\documentclass[18pts]{amsart}
\theoremstyle{plain} \newtheorem{Thm}{ }[section]
\title{Injectivity of the Petri map for twisted Brill-Noether loci}
\author{montserrat teixidor i bigas}
\address{Mathematics Department, Tufts University, Medford MA
02155} \email{montserrat.teixidoribigas@@tufts.edu}

\begin{document}
 \begin{abstract} Let $C$ be a generic curve, $E$ a generic vector bundle on $C$.
  Then, for every line bundle on $C$ the twisted Petri map
  $$P_{E}:H^0(C,L\otimes E)\otimes H^0(C, K\otimes L^*\otimes
E^{*})\rightarrow H^0(C,  K)$$ is injective.

Mathematical sciences classification number 14H60.
 \end{abstract}
\maketitle
\begin{section}{Introduction}
Let $C$ be an algebraic curve defined over an algebraically closed
field. Let $U(r_F,d_F)$ be the moduli space of stable vector bundles
of rank $r$ and degree $d$. Choose a generic vector bundle $E$ of
rank $r_E$ and degree $d_E$ and denote by $B^k_{r_F,d_F}(E)$ the so
called Brill-Noether  locus consisting of vector bundles
$$\{F\in U(r_F,d_F)| \dim H^0(C,F\otimes E)\ge k\}.$$

These sets play a crucial role in the geoemetry of moduli spaces of
vector bundles. For example, for suitable choices of $E$ and $k=1$
they give the generalized theta divisor providing generators of the
Picard group of the moduli space.

Brill-Noether loci can be given scheme structures as locally defined
determinantal varieties (see \cite{GT} for an exposition under the
additional assumption $E={\mathcal O}$). As such their expected
dimension is given by the Brill-Noether number
$$\rho=r_F^2(g-1)+1-k(k-r_Fd_E-r_Ed_F+r_Er_F(g-1))=\dim (U(r_F,d_F)-h^0(F\otimes
E)h^1(F\otimes E).$$ Here we write $h^0(F)$ for the dimension of the
space of sections $H^0(C,F)$ and the last equality in the formula
above assumes that $h^0(F\otimes E)=k$. One also expects that when
$\rho$ is negative, these loci are empty.

Consider the so-called Petri map
$$P_{E}(F):H^0(C,F\otimes E)\otimes H^0(C, K\otimes F^*\otimes
E^{*})\rightarrow H^0(C, F\otimes K\otimes F^*)$$ obtained as the
composition of the natural cup-product and tensorization with the
identity morphism in $E$. If $h^0(F\otimes E)=k$, the tangent space
to $B^k_{r,d}(E)$ at the point $F$ can be identified to the
orthogonal to the image of this map. In particular, $F$ is a
non-singular point of a component of dimension $\rho $ of
$B^k_{r,d}(E)$ if and only if the Petri map is injective. Therefore,
proving the injectivity of the Petri map for the generic curve would
provide a complete positive answer to the more pressing questions in
Brill-Noether Theory. It proves also that when the Brill-Noether
number is negative, the locus is actually empty. Moreover, the
injectivity of the Petri map helps explain the structure of the
tangent cone to the Brill-Noether loci (see \cite{Y}).

In the special case in which $r=1,E={\mathcal O}$ one recovers the
classical Brill-Noether loci in the Jacobian. It is well known that
in this particular case all our expectations are satisfied: the loci
are non-empty on any curve when $\rho >0$ and the Petri map is
injective on the generic curve which implies that these loci have
dimension precisely $\rho$ and that the singular locus of
$B^k_{r,d}(E)$ is $B^{k+1}_{r,d}(E)$ (see for instance \cite{ACGH}).
On the other hand, when $r>1$, there are counterexamples of
particular values of $r,d,k,g$ where the expected results fail, even
for the generic curve (see \cite{GT}).

In this paper, we want to deal with the case $r_F=1$ but arbitrary
$r_E$. The non-emptiness result under these hypothesis was proved by
Ghione in \cite{G}. Therefore, in order to complete the picture, we
need to show the injectivity of the Petri map. We will write $L$
instead of $F$. We can identify $ L\otimes K\otimes L^*$ with $K$
(as it is done for classical Bril-Noether). We prove the following:

\begin{Thm}\label{teorema} {\bf Theorem}
 Given a generic curve
$C$ and a generic  vector bundle $E$ on $C$  for any choice of $L\in
Pic^d(C)$ on $C$, the Petri  map
$$P_E(L):H^0(C,L\otimes E)\otimes H^0(C, K\otimes L^*\otimes
E^{*})\rightarrow H^0(C,K)$$ is injective. In particular, the
twisted Brill-Noether locus $B^k_{1,d}(E)$ is of the expected
dimension $\rho$ and singular only along $B^{k+1}_{1,d}(E)$
\end{Thm}

In order to prove the injectivity of the Petri map  for a generic
curve, it suffices to prove it for a special curve. We'll choose our
curve to be reducible with components  rational and elliptic. Vector
bundles on these curves are easy to describe in terms of the
restrictions of the vector bundles to the various components and
gluing at the nodes. We choose a generic such $E$ and by means of
limit linear series prove the result. The tools used are  those
developed in \cite{petri} which in turn generalizes \cite{EH}.
\end{section}

\begin{section}{Preliminaries on reducible curves}

 Consider a
family of curves $\pi :{\mathcal C}\rightarrow T$. Let $T$ be the
spectrum of a discrete valuation ring ${\mathcal O}$ with maximal
ideal generated by $t$. Assume that the generic fiber of $\pi$ is a
non-singular curve and the special fiber $C$ looks as follows:

Take $g$ elliptic curves $C^i, i=1,\ldots g$ and let $P^i,Q^i$ be
generic points on $E^i$. Take any number of rational curves
$C^0_1,..C^0_{k_0},...C^g_1...C^g_{k_g}$ again with points $P^i_j,
Q^i_j$ on them. Glue $C^i_j$ to $C^i_{j+1}$ by identifying $Q^i_j$
to $P^i_{j+1}$. Glue $C^{i-1}_{k_{i-1}}$ to $C^i$ by identifying
$Q^{i-1}_{k_{i-1}}$ to $P^i$. Glue $C^i$ to $C^i_1$ by identifying
$Q^i$ to $P^i_1$.

 For convenience of notation, we shall
denote by
$$Y_1,...Y_M,\ M=k_0+...+k_g+g$$  the components of $C$ starting with $C^0_1$
and ending with $C^g_{k_g}$. We shall denote by $P_i,Q_i$ the two
points in $Y_i$ that get identified to $Q_{i-1}\in Y_{i-1}$ and
$P_{i+1}\in Y_{i+1}$ respectively. We will use superindices  when we
need to refer to the elliptic components, so $C^i$ will be the
$i^{th}$ elliptic curve..

Note that the form of the central fiber does not change if we make
base changes and normalizations.

Given a vector bundle $E$ and a line bundle $L$ on the generic
fiber, we can assume that we have extensions ${\mathcal E}$,
${\mathcal L}$ to the whole family.

If we tensor ${\mathcal E}\otimes {\mathcal L}$ with a line bundle
of the form ${\mathcal O}_{\mathcal C}(\sum \lambda _iY_i)$, the
restriction of the vector bundle to the central fiber does not
change while the restriction to special fibers changes its degree in
multiples of $r$. In this way,we can concentrate most of the degree
and therefore all of the sections on one component of our choice.
This is the idea behind the definition of limit linear series (see
\cite{EH1}, \cite{petri})

\begin{Thm}\label{lls} {\bf Definition. Limit linear series}
A limit linear series of rank $r$, degree $d$ and dimension $k$ on a
chain of $M$ (not necessarily rational and elliptic) curves consists
of data I,II below for which data III, IV exist satisfying
conditions a-c.

I) For every component $Y_i$, a vector bundle $E_i$ of rank $r$ and
degree $d$ and a $k$-dimensional space of sections $V_i$ of $E_i$.

II) For every node obtained by gluing $Q_i$ and $P_{i+1}$ an
isomorphism of the projectivisation of the fibers $(E_i)_{Q_i}$ and
$(E_{i+1})_{P_{i+1}}$

III) A positive integer a

IV) For every node obtained by gluing $Q_i$ and $P_{i+1}$, bases
$s^t_{Q_i}, s^t_{P_{i+1}}, t=1...k$ of the vector spaces $V_i$ and
$V_{i+1}$

Subject to the conditions

a) $\sum _{i=1}^M d_i-r(M-1)a=d$

b) The orders of vanishing at $P_{i+1},Q_i$ of the sections of the
chosen basis satisfy $ord_{P_{i+1}}s_{i+1}^t+ord_{Q_{i}}s_{i}^t\ge
a$

c) Sections of the vector bundles $E_i(-aP_i), E_i(-aQ_i)$ are
completely determined by their value at the nodes.
\end{Thm}

Note that up to replacing  $E$ by $E\otimes L_0$ and $L$ by
$L\otimes L_0^{-1}$, we can assume that $L$ has degree zero, so this
will be understood in what follows. We then have to deal with only
one rank and degree, so we will drop the subindex $E$ and write
$r,d$ for the degree and rank of $E$.

Write $d=rd_1+d_2,\ 0\le d_2<r$. Denote by $h$ the greatest common
divisor of $r,d$ and write $ r=hr',\ d_2=hd'_2$.

In order to prove the result for a generic $E$, it suffices to prove
it for a particular $E$. We take as vector bundle on the special
curve one with  restriction to the first $g-1$ elliptic components
and all of the rational components a direct sum of line bundles of
the same degree while the restriction to the last elliptic component
is a direct sum of $h$ stable vector bundle of same  rank and
degree.

Choose a component $Y_i$. Modifying the vector bundle on the central
fiber by tensoring with line bundles of the form ${\mathcal
O}_{\mathcal C}(\sum \lambda _iY_i)$, we obtain a vector bundle
${\mathcal E}_{i}$ such  that
$$Y_i\not=C^g , \ \deg ({\mathcal E}_{i|Y_i})=rd_1, \deg ({\mathcal E}_{i|Y_j})=0, \ Y_j\not= Y_i, \ Y_j \not= C^g ,\
 \deg ({\mathcal E}_{i|C^g})=d_2$$
$$Y_i=C^g , \ \deg ({\mathcal E}_{i|Y_i})=d, \deg ({\mathcal E}_{i|Y_j})=0, \ Y_j\not= Y_i$$

Then, the $E_i$ that appears in II of \ref{lls} is ${\mathcal
E}_{i|Y_i}$. The  $a$ of III in \ref{lls} is $a=d_1$.

We do the same procedure with  $K\otimes E^*\otimes L^*$. In this
case, write ${\mathcal E'}$ for ${\mathcal \omega}\otimes{\mathcal
E}^*\otimes{\mathcal L}^*$. Then,

$$Y_i\not=C^g , \ \deg ({\mathcal E}'_{i|Y_i})=r(2g-2-d_1), \deg ({\mathcal E}'_{i|Y_j})=0, \ Y_j\not= Y_i, \ Y_j \not= C^g ,\
 \deg ({\mathcal E}_{i|C^g})=-d_2$$
$$Y_i=C^g , \ \deg ({\mathcal E}'_{i|Y_i})=r(2g-2)-d, \deg ({\mathcal E}'_{i|Y_j})=0, \ Y_j\not= Y_i$$

Then, $a'=2g-3-d_1$.

\begin{Thm}\label{can} {\bf Lemma} The canonical limit linear
series on $C$ has line bundles on $C^i$ equal to $K^i={\mathcal
O}(2(i-1)P^i+2(g-i)Q^i)$ while on the rational components the line
bundle is ${\mathcal O}(2(g-1))$. The space of sections on $C^i$ is
$H^0(L_i(-(i-2)P-(2g-2i)Q))\oplus H^0(L_i(-(2i-1)P-(g-i-1)Q))$. The
unique section whose order of vanishing at $P$ and Q is $2g-2$
  vanishes with order $2(i-1)$ at $P$ and $2g-2i$ at $Q$.
  \end{Thm}
\end{section}

\begin{section}{Vanishing at the nodes of elements in the kernel}

 Let us assume that the twisted Petri map is not injective on the generic curve.
 Then we can find  sections in the kernel of the Petri maps
 $$\pi_*({\mathcal E}\otimes{\mathcal L})_i\otimes \pi_*(\omega _{\pi}\otimes {\mathcal E}^*\otimes{\mathcal
 L}^*)_i\rightarrow \pi_*{\omega_{\pi}}_i$$
 that glue with each other as we go from one component to the next.
 We next study the vanishing of these sections at the various nodes.

\begin{Thm}
{\bf Lemma}.   (see \cite{EH} 1.2 and \cite{petri} 4.1). For every
component $Y_i$, there is a basis $\sigma_j, j=1...k$ of
$\pi_*{\mathcal E}_i$ such that

a) $ord_{P_i}(\sigma_j)=a_j(P_i)$

b) for suitable integers $\alpha_j$, $t^{\alpha_j}\sigma _j$ are a
basis of $\pi _*({\mathcal E}_{i+1})$
\end{Thm}

\begin{Thm}\label{anul}
{\bf Proposition} (see\cite{petri} 4.2) Consider a vector bundle
${\mathcal F}$ on $\pi :{\mathcal C}\rightarrow T$ as before
associated to a vector bundles $F$. Let $\sigma_j$ be a basis of
$\pi_*{\mathcal F}_{i}$ such that $t^{\alpha _j}\sigma _j$ is a
basis of $\pi_*({\mathcal F}_{i+1})$. Let $\bar a$ be as in II of
\ref{lls}. Then , the orders of vanishing of the $\sigma_j$ at the
nodes satisfy

a) $ord_{P_i}(\sigma_j)\le \bar a-ord_{Q_i}\sigma_j\le \alpha _j\le
ord_{P_{i+1}}t^{\alpha_j}\sigma_j$ if $F_i$ is a direct sum of $h$
indecomposable vector bundle of rank $r'$ and  degree $r'\bar a+\bar
d_2, 0< \bar d_2<r'$.

b) $ord_{P_i}(\sigma_j)\le \bar a-ord_{Q_i}\sigma_j\le \alpha
_j-1\le ord_{P_{i+1}}t^{\alpha_j}\sigma_j-1$ if $F_i$ is a direct
sum of $h$ indecomposable vector bundle of rank $r'$ and  degree
$r'(\bar a-1)+\bar d_2, 0<\bar d_2<r'$.

c) $ord_{P_i}(\sigma_j)\le \bar a-ord_{Q_i}\sigma_j\le \alpha _j\le
ord_{P_{i+1}}t^{\alpha_j}\sigma_j$ if $F_i$ is  a direct sum of line
bundles of degree $\bar a$.

Moreover, if equality holds, then $\sigma_j$  as a section of $F_i$
  vanishes at
 $P_i,Q_i$ with orders adding up to $\bar a$ in a), c) and $\bar a -1$ in b).

\end{Thm}

Consider now the Petri map
$$\pi_*({\mathcal L}\otimes {\mathcal E})\otimes \pi_*(\omega_{\pi}\otimes {\mathcal L}^{-1}\otimes {\mathcal
E}^{*}) \rightarrow \pi_* {\omega _{\pi}}$$

As in \cite{EH}, p.277, one can define the order of  a section $\rho
\in \pi_*{\mathcal L}\otimes {\mathcal E}_Y\otimes \pi_*(\omega
_{\pi}\otimes {\mathcal L}^{-1}\otimes {\mathcal E}^{*})$ at a point
$P$ on a component $Y$ as follows:

\begin{Thm} {\bf Definition} We say $ord_P(\rho _{|Y})\ge l$ if and only if $\rho$ is in the
linear span of $t(\pi_*{\mathcal L}\otimes {\mathcal E}_Y\otimes
\pi_*(\omega_{\pi}\otimes {\mathcal L}^{-1}\otimes {\mathcal
E}^{*}))$ and elements of the form $\sigma_m\otimes \sigma'_n$ where
$ord_P(\sigma_m)+ord_P(\sigma' _n)\ge l,\ \sigma_m \in
\pi_*({\mathcal L}\otimes {\mathcal E}_Y),\ \sigma '_n \in
\pi_*(\omega_{\pi}\otimes {\mathcal L}^{-1}\otimes {\mathcal
E}^{*})$.
\end{Thm}

One then has the following result (cf. \cite{EH}, Lemma 3.2)

\begin{Thm}\label{anulrho} {\bf Lemma} Let $\sigma_m$ be a basis
of the free ${\mathcal O}$ module $\pi_*({\mathcal L}\otimes
{\mathcal E})_i$ such that the orders of vanishing of the $\sigma
_m$ at $P_i$ are the orders of vanishing of the linear series at
this point and $t^{\alpha _m}\sigma _m$ is a basis of
$\pi_*({\mathcal L}\otimes {\mathcal E})_{i+1}$. Let $\sigma'_n$ be
a basis of the free ${\mathcal O}$ module $\pi_*(\omega_{\pi}\otimes
{\mathcal L}^{-1}\otimes {\mathcal E}^{*}))$ such that the orders of
vanishing of the $\sigma' _n$ at $P_i$ are the distinct orders of
vanishing of the linear series at this point and $t^{\alpha '
_n}\sigma' _n$ is a basis of $\pi_*((\omega_{\pi}\otimes {\mathcal
L}^{-1}\otimes {\mathcal E}^{*}))_{i+1})$. If
$$\rho =\sum f_{n,m}(\sigma_m\otimes\sigma'_n)$$ where the $f_{n,m}$ are functions on the discrete valuation
ring ${\mathcal O}$ and the associated discrete valuation is $\nu$,
then
$$ord_{P_i}(\rho_{|Y_i})=min_{\{\nu(f_{n,m})=0\} } (ord_{P_i}(\sigma_m)+ord_{P_i}(\sigma'_n))$$
If $\beta_i$ is the unique integer such that
$$t^{\beta _i}\rho \in \pi _*({\mathcal L}_i\otimes {\mathcal E})_i\otimes
 \pi_*(\omega_{\pi}\otimes {\mathcal L}^{-1}\otimes
 {\mathcal E}^{*})_{i}-t(\pi _*{\mathcal L}\otimes {\mathcal E}_i\otimes
\pi_*(\omega_{\pi}\otimes {\mathcal L}^{-1}\otimes {\mathcal
E}^{*})_{i})$$ then
$$\beta_i =max\{ \alpha _m+\alpha ' _n-\nu (f_{nm})\}$$
\end{Thm}

Let us assume now that the kernel of the Petri map  is non-zero on
the generic curve. We can then find an element $\rho$ such that say
$$t^{\beta _i}\rho \in \pi_*({\mathcal L}\otimes {\mathcal E})_i\otimes
\pi_*(\omega_{\pi}\otimes {\mathcal L}^{-1}\otimes {\mathcal
E}^{*})_{i}-t(\pi _*({\mathcal L}\otimes {\mathcal E})_i\otimes
\pi_*(\omega_{\pi}\otimes {\mathcal L}^{-1}\otimes {\mathcal
E}^{*})_{i})$$ and $t^{\beta _{i+1}}\rho \in  \pi_*({\mathcal
L}\otimes {\mathcal E})_{i+1}\otimes \pi_*(\omega_{\pi}\otimes
{\mathcal L}^{-1}\otimes {\mathcal E}^{*})_{i+1}$ and in the kernel
of the twisted Petri map to $ \pi_*(\omega_{\pi})_{i+1}$.

\begin{Thm} {\bf Proposition} \label{claim}

1) If $Y_i$ is any component,
$$ord_{P_{i+1}}(t^{\beta_{i+1}}\rho_{|Y_{i+1}})\ge
ord_{P_i}(t^{\beta_i}\rho_{|{Y_i}})$$.

2) If $Y_i$ is an elliptic component, $Y_i\not= C^g$
$$ord_{P_{i+1}}(t^{\beta_{i+1}}\rho_{|Y_{i+1}})\ge
ord_{P_i}(t^{\beta_i}\rho_{|{Y_i}})+1$$

If $Y_i$ is the first component for which there is equality in the
above inequality, then the terms of $\rho_i$ that give the vanishing
at $Q_i$ can be written in the form
$$(*)\bar\sigma\otimes \sigma'+\sigma\otimes\bar\sigma'$$ for two
specific sections  $\bar \sigma \in({\mathcal L}\otimes {\mathcal
E})_i$ and $ \bar \sigma' \in\pi_*(\omega_{\pi}\otimes {\mathcal
L}^{-1}\otimes {\mathcal E}^{*})_{i}$. Moreover, the fibers of
$\sigma, \sigma '$ at the nodes move in spaces of dimension at most
$r-1$ complementary to the directions of $\bar \sigma,\
\bar\sigma'$.

If $\rho_{i-1}$ is of the form (*), then
$$ord_{P_{i+1}}(t^{\beta_{i+1}}\rho_{|Y_{i+1}})\ge
ord_{P_i}(t^{\beta_i}\rho_{|{Y_i}})+2$$ and if there is equality in
this inequality, then (*) holds again in the component $Y_i$..

3) If $Y_i$ is the last elliptic component, $Y_i= C^g$
$$ord_{P_{i+1}}(t^{\beta_{i+1}}\rho_{|Y_{i+1}})\ge
ord_{P_i}(t^{\beta_i}\rho_{|{Y_i}})+1$$  If $\rho_{i-1}$ is of the
form (*), then the inequality above is strict.
\end{Thm}

\begin{proof}

 We  first prove that  1) holds on any component.  Choose a
basis $\sigma _m,\ m=1...k$ of $\pi_*({\mathcal L}\otimes {\mathcal
E}_i)$ such that $t^{\alpha _m}\sigma _m$ is a basis of
$\pi_*({\mathcal L}\otimes {\mathcal E}) _{i+1}$.  Similarly, choose
a basis $\sigma ' _n,\ n=1...k'=k-rd+r(g-1)$ of
$\pi_*(\omega_{\pi}\otimes{\mathcal L}^{-1}\otimes {\mathcal E})_i$
such that $t^{\alpha' _n}\sigma' _n$ is a basis of
$\pi_*(\omega_{\pi}\otimes {\mathcal L}^{-1}\otimes {\mathcal E}^{*}
_{i+1})$. For simplicity of notation, we shall assume that
$\beta_i=0$. Write
$$\rho =\sum _{ m,n}f_{m,n}(\sigma_m\otimes \sigma '_n) $$ Then, from
\ref{anulrho},
$$ord_{P_i}(\rho _{|Y_i})=min_{\{ \nu
(f_{m,n})=0\}}(ord_{P_i}(\sigma_m)+ord_{P_i}(\sigma'_n))$$ Assume
that this minimum is attained by a pair corresponding to the indices
$m_0,n_0$ with $\nu(f_{m_0,n_0})=0$. Then from \ref{anul},
$$(ord_{P_i}(\sigma_{m_0})+ord_{P_i}(\sigma'_{n_0}))\le
a+a'- ord_{Q_i}(\sigma_{m_0})-ord_{Q_i}(\sigma'_{n_0}) \le
\alpha_{m_0}+\alpha'_{n_0}$$ From \ref{anulrho} and the fact that
$\nu (f_{n_0,m_0})=0$, the latter is at most $\beta_{i+1}$.

Write $$t^{\beta_{i+1}}\rho =\sum _{n\le
m}(t^{\beta_{i+1}-\alpha_m-\alpha'_n}f_{nm})(t^{\alpha_m}\sigma_m\otimes
t^{\alpha'_n}\sigma '_n )$$

 Hence, from \ref{anulrho} $$ord_{P_{i+1}}(t^{\beta_{i+1}}\rho
_{|Y_{i+1}})=min_{\{ \beta_{i+1}-\alpha_m-\alpha'_n+\nu
(f_{nm})=0\}}(ord_{P_{i+1}}(t^{\alpha_m}\sigma_m)+ord_{P_{i+1}}(t^{\alpha
_n}\sigma_n))$$

Assume that this minimum is attained at a pair $m_1,n_1$ with
$$\beta_{i+1}-\alpha_{m_1}-\alpha'_{n_1}+\nu(f_{m_1,n_1})=0$$
Then, $$\beta_{i+1}\le
\beta_{i+1}+\nu(f_{m_1,n_1})=\alpha_{m_1}+\alpha'_{n_1}\le
ord_{P_{i+1}}(t^{\alpha_{m_1}}\sigma_{m_1})+ord_{P_{i+1}}(t^{\alpha'
_{n_1}}\sigma'_{n_1}))$$ where the last inequality comes from
\ref{anul}

 Stringing together the above inequalities, we obtain
$$ord_{P_i}(\rho_{|Y_i})\le ord_{P_{i+1}}(\rho_{|Y_{i+1}}).$$ Hence
part 1) is proved.

 Assume now that there is equality in the
inequality above and we are in the situation of b). Then all the
previous inequalities must be equalities. In particular, any terms
$\sigma_{m}\otimes \sigma'_{n}$ that give the vanishing of $\rho$ at
$P_i$ satisfy
$$ord_{P_i}(\sigma _m)+ord_{Q_i}(\sigma_m)=a,\ ord_{P_i}(\sigma' _n)+ord_{Q_i}(\sigma'_n)=a'.$$

Write ${\mathcal E}_i$ as a direct sum of $r$ line bundles of degree
$d_1$
$${\mathcal E}_i=L^i_1\oplus \ldots \oplus L^i_r$$
By the genericity of $E'$, the $L^i_j$ are generic. Then,
$$({\mathcal E}\otimes {\mathcal L})_i=(L^i_1\otimes L^i)\oplus \ldots \oplus (L^i_r\otimes L^i).$$
and at most one of the $L^i_j\otimes L^i$ is of the form ${\mathcal
O}(a^i_jP_i+(d_1-a^i_j)Q_i)$. Denote this index (if it exists) by
$j_i$. On an elliptic curve, with $P_i,\ Q_i$ generic points, there
is only one section $\bar \sigma^i$ of a line bundle of degree $d_1$
with orders of vanishing at $P_i,\ Q_i$ adding up to $d_1$, the one
corresponding to $a^i_jP_i+(d_1-a^i_j)Q_i$. Hence, there is only one
section of $({\mathcal E}\otimes {\mathcal L})_i$ with orders of
vanishing at $P_i,\ Q_i$ adding up to $d_1$, the one corresponding
to $a^i_jP_i+(d_1-a^i_j)Q_i$ on the $j_i^{th}$ component and zeros
elsewhere.

From the description of the limit canonical series, at most one of
the $K\otimes L_j^{i*}\otimes L^{i*}$ is of  the form ${\mathcal
O}(a'^i_jP_i+(a'-a'^i_j)Q_i)$ and it corresponds to the same $j_i$
as before. Therefore, there is only one section $\bar \sigma'^i$ of
$(\omega_{\pi}\otimes {\mathcal E}^*\otimes {\mathcal L}^*)_i$ with
orders of vanishing at $P_i,\ Q_i$ adding up to $2g-2-d_1$ and again
the only non-zero part of this section corresponds to the summand
$j_i$ in the direct sum decomposition of the vector bundle.

  As $\bar \sigma^i\otimes \bar \sigma'^i$ is not in
the kernel of the Petri map, the vanishing must go up by at least
one and if it goes up by exactly one, the section $\rho$ is of the
form given in (*). Moreover, the $j_i^{th}$ component of $\bar
\sigma^i ,\bar \sigma'^i$ is a section of the canonical bundle that
vanishes only at $P_i,Q_i$. It follows from \ref{can} that it
vanishes at $P_i$ to order $2(i-1)$.

Assume now that this $Y_i$ is the first elliptic component $C^m$
where the vanishing goes up by just one. Then,  $$ord_{P_i}(\rho)\ge
2m-2.$$ Write $\sigma=(\sigma_1,\ldots, \sigma_r), \sigma
'=(\sigma'_1,\ldots \sigma'_r)$ according to the decomposition of
${\mathcal E}_i$. We want to prove that if the vanishing goes up by
exactly one, then $\sigma_{j_i}=0,\ \sigma'_{j_i}=0$. This will show
that the fibers of $\sigma, \sigma '$ at the nodes move in  spaces
of dimension at most $r-1$ complementary to the directions of $\bar
\sigma,\ \bar \sigma'$.

By assumption, $ ord_{P_i}(\rho)=ord_{P_i}(\bar
\sigma^i)+ord_{P_i}(\sigma')= ord_{P_i}(\sigma)+ord_{P_i}(\bar
\sigma'^i)$. Assume that the $j_i$ component of $\sigma'$ is
non-zero. Then, both this and the $j_i$ component of $\bar \sigma'$
 are sections of $(L^i_{j_i}\otimes L^i)'$. Hence, they vanish at
 different orders at $P_i$.

 Then
$$2m-2\le ord_{P_i}(\rho)=ord_{P_i}(\bar \sigma^i)+ord_{P_i}(\sigma') \not=
ord_{P_i}(\bar \sigma^i)+ord_{P_i}(\bar \sigma'^i)=2m-2.$$ Hence,
$ord_{P_i}(\bar \sigma^i)+ord_{P_i}(\sigma')\ge 2m-1$.

From \ref{anul}, $$ord_{P_i}(\sigma ')\le 2g-2-d-ord_{Q_i}(\sigma '
)\le \alpha '\le ord_{P_{i+1}}(t^{\alpha}\sigma ').$$

The first inequality is strict. Hence,
$$ord_{P_{i+1}}(t^{\alpha _i}\rho)=ord_{P_{i+1}}(t^{\bar \alpha}
\bar\sigma^i)+ord_{P_{i+1}}(t^{\alpha' }\sigma ')\ge$$
$$ord_{P_i}(\bar\sigma^i)+ord_{P_i}(\sigma ')+1\ge 2m-1+1=2m.$$

If $\rho_{i-1}$ is of the form (*), the genericity of the gluing at
the node between the components $Y_{i-1}$ and $Y_i$ means that two
sections of the form (*) corresponding to the sections with highest
vanishing in each component will not glue together. Therefore, the
vanishing in the next component increases in at least two.

Assume now that it increases in precisely two units. By the
genericity of the gluing, none of the the sections that glue with
each of $\bar\sigma^{i-1},\ \sigma'_{i-1},\ \sigma^{i-1},
\bar\sigma'^{i-1}$ can be $\bar\sigma^{i}, \bar\sigma^{'i}$.
Therefore, for each of them the sum of vanishing at the two nodes
must be one less than the maximum, namely $a-1$. As $E^i$ is a
direct sum of line bundles of degree $a$, there is a unique section
up to scalars whose order of vanishing at the two nodes is $a-1$ and
glues with a preassigned direction at $P_i$. Hence, the fiber at
$Q_i$ of $\rho$ is completely determined and the last statement of
2) is proved.

Consider now the last elliptic component. Then, the sum of the
orders of vanishing between $P_i,Q_i$ of a section of  ${\mathcal
E}_i$ is at most $a$ and there is a space of dimension $d_2$ of such
sections. The sum of the orders of vanishing between $P_i,Q_i$ of a
section of  ${\mathcal E'}_i$ is at most $a'-1$ and there is a space
of dimension $r-d_2$ of such sections. From the arguments above,
this implies that the vanishing increases in at least one. Moreover,
if it increases in exactly one, then $\rho_i$ can be written in the
form $\sum \sigma_k\otimes \sigma'_k$ where all the $\sigma_k,\
\sigma'_k$ are in the special directions mentioned above. In
particular, this $\rho$ does not glue with one of the form in (*).
Hence, 3) is proved.
\end{proof}

The proof of the Theorem now follows from the previous proposition :
we have shown that on rational components the vanishing at the nodes
does not decrease while on elliptic components it increases in at
least one and in all but one component it increases in at least two.
This implies that the vanishing at one point in the last component
is at least $2g-1$. This is impossible, as the canonical line bundle
has degree $2g-2$.

\end{section}

\end{document}